\documentclass[conference]{IEEEtran}
\IEEEoverridecommandlockouts
\ifCLASSINFOpdf
  \usepackage[pdftex]{graphicx}
  \graphicspath{{./Images/}}
  \DeclareGraphicsExtensions{.pdf,.jpeg,.png}
\else
  \usepackage[dvips]{graphicx}
  \graphicspath{{./Images/}}
  \DeclareGraphicsExtensions{.eps}
\fi
\usepackage{color}
\usepackage{amsmath}
\usepackage{amsfonts}
\usepackage{amsthm}
\usepackage{amssymb}
\usepackage{multirow}

\newtheorem{theorem}{Theorem}

\usepackage{algorithm}
\usepackage{algorithmic}

\DeclareMathOperator{\diag}{diag}
\DeclareMathOperator{\prox}{prox}

\def\st{\mbox{s.t.}}
\newcommand{\N}{\mathbb{N}}
\newcommand{\R}{\mathbb{R}}
\newcommand{\C}{\mathbb{C}}
\def\bb{\mathbf b}

\def\bd{\mathbf d}

\def\bg{\mathbf g}
\def\bp{\mathbf p}
\def\bq{\mathbf q}
\def\bs{\mathbf s}
\def\bu{\mathbf u}
\def\bv{\mathbf v}
\def\bw{\mathbf w}
\def\bx{\mathbf x}
\def\by{\mathbf y}
\def\bz{\mathbf z}
\def\b0{\mathbf 0}

\def\bgamma{\boldsymbol \gamma}
\def\bmu{\boldsymbol \mu}
\def\bxi{\boldsymbol \xi}

\def\mE{\mathcal{E}}
\def\mF{\mathcal{F}}
\def\mL{\mathcal{L}}

\begin{document}
%
\title{TGV-based restoration of Poissonian images with
automatic estimation of the regularization parameter
\thanks{This research was supported by the Istituto Nazionale di Alta Mate\-ma\-tica, Gruppo Nazionale per il Calcolo Scientifico (INdAM-GNCS). DdS and MV were also funded by the V:ALERE Program of the University of Campania ``L.~Vanvitelli''.}}

\author{\IEEEauthorblockN{Daniela di Serafino}
\IEEEauthorblockA{Department of Mathematics and Applications\\
University of Naples Federico II\\
Naples, Italy\\
Email: daniela.diserafino@unina.it}
\and
\IEEEauthorblockN{Germana Landi}
\IEEEauthorblockA{Department of Mathematics\\
University of Bologna\\
Bologna, Italy\\
Email: germana.landi@unibo.it}
\and
\IEEEauthorblockN{Marco Viola}
\IEEEauthorblockA{Department of Mathematics and Physics\\
University of Campania ``L. Vanvitelli''\\
Caserta, Italy\\
Email: marco.viola@unicampania.it}}


\IEEEspecialpapernotice{VERSION 1 -- April 25, 2021}

\maketitle

\begin{abstract}
The problem of restoring images corrupted by Poisson noise is common in many application fields and, because of its intrinsic ill posedness, it requires regularization techniques for its solution. The effectiveness of such techniques depends on the value of the regularization parameter balancing data fidelity and regularity of the solution. Here we consider the Total Generalized Variation regularization introduced in [\emph{SIAM J. Imag. Sci, 3(3), 492–526, 2010}], which has demonstrated its ability of preserving sharp features as well as smooth transition variations, and introduce an automatic strategy for defining the value of the regularization parameter. We solve the corresponding optimization problem by using a 3-block version of ADMM. Preliminary numerical experiments support the proposed approach.
\end{abstract}


%
\IEEEpeerreviewmaketitle

\section{Introduction}
We are interested in the problem of restoring images corrupted by Poisson noise, which arises in many application areas, e.g., in fluorescence microscopy~\cite{sarder:2006}, computed tomography (CT) \cite{herman:2009}, and astronomical imaging~\cite{bertero:2009}. This inverse problem is usually highly ill conditioned and regularization techniques are required in order to obtain reasonable approximate solutions. Such techniques reformulate the image restoration problem as a minimizzation problem whose objective function contains a data fidelity term and a regularization term, imposing some a-priori information on the object to be restored. These two terms are balanced by a regularization parameter whose value greatly affects the quality of the restored image. Additional constraints can be incorporated in the problem to better model the image. 
When the noise in the data has a Poisson distribution, the Kullback-Leibler divergence is usually considered as data fidelity term and nonnegativity constraints are imposed.
Several regularization terms have been proposed in the literature. For real-life images one of the most successful regularizers is the Total Generalized Variation (TGV), introduced in \cite{bredies:2010} as a way to overcome the well-known staircasing effect of the classical Total Variation (TV) \cite{rudin:1992}. TGV can be seen as a higher-order counterpart of TV. Indeed, its aim is to enforce higher-order smoothness in the restored image while maintaining the ability of preserving sharp edges. 
However, the selection of an appropriate value for the regularization parameter is a crucial issue. Even if several criteria for choosing the parameter have been discussed in the literature when the noise in the image follows a Gaussian distribution, to the best of our knowledge only a few methods have been proposed to select the value of the regularization parameter for Poisson noise. In \cite{bertero:2010}, a discrepancy principle for Poisson noise was proposed and theoretically analysed. In \cite{Bardsley_2009} a similar principle was derived by using a quadratic approximation of the data fidelity function. The generalized cross validation criterion and the unbiased predictive risk method were also introduced in \cite{Bardsley_2009} by using the same quadratic approximation. A modification of the discrepancy principle of \cite{bertero:2010} and \cite{Bardsley_2009} was discussed in \cite{carlavan2012}. In \cite{Chen_2011} a spatially-adapted regularization parameter selection scheme was derived by using a confidence interval technique based on the expected maximal local variance estimate. 

In this work, inspired by \cite{bortolotti:2016} and \cite{ito:2011}, we present an automatic strategy for the determination of the regularization parameter balancing the two terms of the objective function. Then we discuss a 3-block version of the Alternating Direction Method of Multipliers (ADMM) \cite{boyd:2011admm, cai:2017} for the solution of the nonsmooth bound-constrained optimization problem modeling the image restoration. Finally, we show the results of preliminary numerical experiments using the proposed approach on both synthetic and real-life images.


\section{The KL-TGV$^2$ model}

A discrete formulation of the Poisson restoration problem is given next. Let the image to be restored be represented by $\bu \in \R^n$, where $n = n_1 n_2$ and $n_1 \times n_2$ is the image size in pixels, and let the measured data be a vector $\bb \in \N_0^n$. Since the pixels in the image are related to the number of photons hitting sensors, the entries of $\bb$ can be considered as samples from $n$ independent Poisson random variables 
$B_j$ with probability
$$
    P(B_j = b_j) = \frac{e^{-[A \bu + \bgamma]_j}[A\bu + \bgamma]_{j}^{b_j}}{b_j!}.
$$
Here the matrix $A = (a_{ij}) \in \R^{n \times n}$ models the blur affecting the image and $\bgamma \in \R^{n}$,  $\bgamma > 0$, models the background radiation detected by the sensors. $A$ is usually required to satisfy the conditions
$$  a_{ij} \geq 0 \mbox{ for all }  i, j, \qquad \sum_{i=1}^n a_{ij} = 1 \mbox{ for all } j.$$

\noindent
By following the maximum-likelihood approach \cite{bertero:2009}, we can estimate $\bu$
by minimizing the negative logarithm of the likelihood function, i.e.,
by solving the problem
\begin{equation} \label{eq:kl}
    \!\!\!\!\! \min_{\bu\geq 0} D_{KL}(\bu) 
    \equiv \sum_{i=1}^n \left( b_i \ln \frac{b_i}{[A \bu+\bgamma]_i} + [A \bu+\bgamma]_i - b_i \right),
\end{equation}
\noindent
where $D_{KL}(\bu)$ denotes the Kullback-Leibler (KL) divergence of $A \bu + \bgamma$ from $\bb$ and we set
$$ b_i \ln \frac{b_i}{[A \bu+\bgamma]_i} = 0 \quad \mbox{if} \;\; b_i = 0. $$
\noindent
We observe that $D_{KL}$ is a convex function and it is strictly convex only if the kernel of $A$ is trivial.

Among the regularization terms falling in the category of TGV regularization, we consider the second-order regularization TGV$^2$, which has the form
\begin{equation}\label{eq:TGV2}
\begin{split}
    \mathrm{TGV}^2(\bu) = \min\limits_{\bw\in\R^{2n}} \displaystyle 
    {\alpha_0 \left\|\nabla \bu - \bw\right\|_{2,1|\R^{2n}}} \\
    + {\alpha_1 \left\|\mE\bw\right\|_{2,1|\R^{4n}}},
\end{split}
\end{equation}
where $\alpha_0,\alpha_1\in(0,\,1)$ are positive parameters balancing the two components of the regularization. The linear operators $\nabla\in\R^{2n\times n}$ and $\mE\in\R^{4n\times2n}$ represent the \textit{discrete gradient operator} and the \textit{symmetrized derivative operator}, respectively, which are described below. Note that for any $\bv \in \R^{2n}$ and $\by\in\R^{4n}$ we define
\begin{equation}\label{eq:def_21norm_r2n}
  \|\bv\|_{2,1|\R^{2n}} = \sum_{j=1}^n \sqrt{v_j^2+v_{n+j}^2},
\end{equation}
and
\begin{equation}\label{eq:def_21norm_r4n}
  \|\by\|_{2,1|\R^{4n}} = \sum_{j=1}^n \sqrt{y_j^2 + y_{n+j}^2 + y_{2n+j}^2 + y_{3n+j}^2}.  
\end{equation}

Let $D_H, D_V \in\R^{n\times n}$ represent the forward finite-difference operators along the horizontal and the vertical direction, respectively. The discrete gradient operator and the discrete symmetrized derivative operator are defined as
\begin{equation*}
    \nabla = \left[\begin{array}{c}
         D_H \\ D_H \end{array}\right], \qquad 
    \mE = \left[\begin{array}{cc}
         D_H & 0 \\
         \frac{1}{2}D_V & \frac{1}{2}D_H\\
         \frac{1}{2}D_V & \frac{1}{2}D_H\\
         0 & D_V
         \end{array}\right].
\end{equation*}

The combination of \eqref{eq:kl} and \eqref{eq:TGV2} results in the nonsmooth optimization problem
\begin{equation}\label{eq:kltgv2}
    \begin{array}{ll}
    \min\limits_{\bu,\bw} & \displaystyle \lambda\,D_{KL}(\bu)  +  \alpha_0 \left\|\nabla \bu - \bw\right\|_{2,1|\R^{2n}} + \alpha_1 \left\|\mE\bw\right\|_{2,1|\R^{4n}} \\
    \st       & \bu \geq 0,
    \end{array}
\end{equation}
\noindent
which we refer to as KL-TGV$^2$ model, where $\lambda\in(0,\,+\infty)$ is the parameter regulating the relative weight between the KL divergence and the two terms composing the TGV$^2$ regularization.

\section{An automatic strategy for the computation of the regularization parameter}
Recently, two closely related principles have been proposed for selecting the values of the regularization parameters in multipenalty regularization without any a-priori information on the noise norm. The \emph{uniform penalty principle} \cite{bortolotti:2016} selects the regularization parameters to ensure uniform penalization, i.e., a constant value of the all the penalty terms equal to the data fidelity term. This principle has been successfully used for the inversion of Nuclear Magnetic Resonance data by using spatially-adapted $\ell_2$ and $\ell_1-\ell_2$ regularization \cite{bortolotti:2021}. The \emph{balancing principle} \cite{ito:2011} chooses the values of the regularization parameters by balancing the data fidelity and regularization terms; a deep theoretical analysis of this principle is presented in \cite{ito:2014} together with a convergent fixed-point iterative scheme for its realization.

When only one regularization term is considered, both the principles suggest to compute $\bu$ and $\lambda$ such that
\begin{equation}\label{principle}
  \left\{
    \begin{array}{l}
      \displaystyle \bu = \arg\min\limits_{\bu\geq0} \lambda\,D_{KL}(\bu) + {TGV}^2(\bu)  \\
      \displaystyle \lambda = \gamma \frac{{TGV}^2(\bu)}{D_{KL}(\bu)}, 
    \end{array}
  \right.
\end{equation}
where the multiplicative constant $\gamma$ differs in the two principles. 

For computing an approximate solution to \eqref{principle}, the fixed-point type approach sketched in Algorithm~\ref{alg:fixedpoint} can be employed \cite{bortolotti:2016,ito:2014}.
\begin{algorithm}[h]
    \caption{\label{alg:fixedpoint}}
    {\small
    \begin{algorithmic}[1]
        \STATE Choose $\gamma>0$, $\lambda^0$, $\bu^0$, and set $j=0$.
        \REPEAT
        \STATE {$ \displaystyle \bu^{j+1} = \arg\min\limits_{\bu\geq0} \lambda^j\,D_{KL}(\bu) + {TGV}^2(\bu)$ \label{alg:fixedpoint_line_3}}
        \STATE $ \lambda^{j+1} = \gamma \frac{{TGV}^2(\bu^{j+1})}{D_{KL}(\bu^{j+1})}$
        \STATE $ j = j+1 $
        \UNTIL{a stopping criterion is satisfied}
    \end{algorithmic}
    }
\end{algorithm}

\noindent
Note that Algorithm~\ref{alg:fixedpoint} requires a method for solving the KL-TGV$^2$ problem. A method performing this task is described in the next section.

\section{3-block ADMM for the minimization of the KL-TGV$^2$ model \label{sec:3block_admm}}

Now we focus on the solution of problem \eqref{eq:kltgv2}. Henceforth,
we assume that periodic boundary conditions are considered for both the linear operator $A$ and the finite-difference operators $D_H$ and $D_V$. This implies that they have a Block Circulant with Circulant Blocks (BCCB) structure.

Problem \eqref{eq:kltgv2} is a nonsmooth bound-constrained minimization problem, hence its solution cannot be computed by ``standard'' optimization techniques. A possibility for dealing with nonsmoothness is to replace the norms \eqref{eq:def_21norm_r2n} and \eqref{eq:def_21norm_r4n} with smooth approximations and then to solve the resulting problem with first- or second-order smooth optimization methods (see, e.g., \cite{bonettini:2009, birgin:2014, diserafino:2018amc, diserafino:2020amc} and the references therein). It is also possible to address the nonsmooth problem as it is by applying, e.g., methods based on Bregman iterations or splitting techniques (see, e.g., \cite{bregman:1967,goldstein:2009,desimone:2020} and the references therein).

Here we propose to solve problem \eqref{eq:kltgv2} by a 3-block version of ADMM \cite{boyd:2011admm}. First of all, let us reformulate problem \eqref{eq:kltgv2} as follows:
\begin{equation}\label{eq:kltgv2_reformulate}
    \begin{array}{ll}
    \min\limits_{\bu,\bw,\bz_1,\bz_2} & \displaystyle \lambda\,D_{KL}(\bu) + \alpha_0\,\|\bz_0\|_{2,1|\R^{2n}} + \alpha_1\,\|\bz_1\|_{2,1|\R^{4n}}\\
    \st     & \nabla\bu - \bw = \bz_0,\\
            & \mE\bw = \bz_1,\\
            & \bu\geq 0,
    \end{array}
\end{equation}
where we introduce the auxiliary variables $\bz_0\in\R^{2n}$ and $\bz_1\in\R^{4n}$ to separate the smooth part from the nonsmooth ones. By setting $\bz = [\bz_0^\top,\; \bz_1^\top]^\top$, we can write the KL-TGV$^2$ problem~as
\begin{equation}\label{eq:kltgv2_reformulate_admm}
    \begin{array}{ll}
    \min\limits_{\bu,\bw,\bz} & \displaystyle F_u(\bu) + F_w(\bw) + F_z(\bz) \\
    \st     & M_u\,\bu + M_w\,\bw + M_z\,\bz = 0,
    \end{array}
\end{equation}
where
\begin{equation}\label{eq:define_F}
\begin{split}
F_u(\bu)= \lambda\,D_{KL}(\bu)+ \chi_{\R_+^n}(\bu),\quad F_w(\bw) =0,\\
F_z(\bz) = \alpha_0\,\|\bz_0\|_{2,1|\R^{2n}} + \alpha_1\,\|\bz_1\|_{2,1|\R^{4n}},
\end{split}
\end{equation}
and the matrices $M_u\in\R^{6n\times n}$, $M_w\in\R^{6n\times 2n}$ and $M_z\in\R^{6n\times 6n}$ are defined as
\begin{equation}\label{eq:define_M}
M_u = \left[ \begin{array}{c} \nabla \\ 0 \end{array} \right], \quad M_w = \left[ \begin{array}{c} -I_{2n} \\ \mE \end{array} \right], \quad
M_z = -I_{6n}.
\end{equation}

The Lagrangian function of problem \eqref{eq:kltgv2_reformulate_admm} is defined as
\begin{equation}\label{eq:kltgv2_lagrangian}
\begin{split}
    \mL(\bu,\bw,\bz,\bxi) =& F_u(\bu) + F_w(\bw) + F_z(\bz) \\
    &+ \bxi^\top\left(M_u\,\bu + M_w\,\bw + M_z\,\bz\right),
\end{split}
\end{equation}
where $\bxi\in\R^{6n}$ is a vector of Lagrange multipliers,
and the augmented Lagrangian function associated with problem \eqref{eq:kltgv2_reformulate_admm} is defined as
\begin{equation}\label{eq:kltgv2_augmented_lagrangian}
\begin{split}
    \mL_A(\bx,\bz,\bxi;\rho) =& F_u(\bu) + F_w(\bw) + F_z(\bz) \\
    &+ \bxi^\top\left(M_u\,\bu + M_w\,\bw + M_z\,\bz\right)\\
    &+ \frac{\rho}{2}\left\|M_u\,\bu + M_w\,\bw + M_z\,\bz\right\|_2^2,
\end{split}
\end{equation}
where $\rho>0$.

A 3-block ADMM for the solution of problem~\eqref{eq:kltgv2_reformulate_admm} can be formulated as follows. Let $\bu^0\in\R^{n}$, $\bw^0\in\R^{2n}$, $\bz^0\in\R^{6n}$, and $\bxi^0\in\R^{6n}$. At each iteration $k>0$ ADMM determines the iterate $( \bu^{k+1}, \bw^{k+1}, \bz^{k+1}, \bxi^{k+1} )$ as
\begin{equation}\label{eq:admm_method}
\begin{split}
    \bu^{k+1} & = \displaystyle \mathrm{arg}\min\limits_{\bu\in\R^{n}}  \mL_A(\bu,\bw^k,\bz^k,\bxi^k;\rho),\\
    \bw^{k+1} & = \displaystyle \mathrm{arg}\min\limits_{\bu\in\R^{2n}} \mL_A(\bu^{k+1},\bw,\bz^k,\bxi^k;\rho),\\
    \bz^{k+1} & = \displaystyle \mathrm{arg}\min\limits_{\bz\in\R^{6n}} \mL_A(\bu^{k+1},\bw^{k+1},\bz,\bxi^k;\rho),\\
    \bxi^{k+1} & = \displaystyle \bxi^k + \rho\left(M_u\,\bu^{k+1} + M_w\,\bw^{k+1} + M_z\,\bz^{k+1}\right).
\end{split}
\end{equation}

\noindent
We note that the functions $F_u(\bu)$, $F_w(\bw)$ and $F_z(\bz)$ in \eqref{eq:kltgv2_reformulate_admm} are closed, proper and convex; moreover, the matrices $M_w$ and $M_z$ defined in \eqref{eq:define_M} have full (column) rank.

Despite its good practical performance in many application areas, the convergence of the 3-block ADMM scheme in the general case is still an open problem. The main issue in \eqref{eq:kltgv2_reformulate_admm} is that the objective function $F_u(\bu)$ is not strongly convex and the matrix $M_u$ is rank deficient. This implies that the objective function of the first subproblem 
in \eqref{eq:admm_method} is not strongly convex. This issue can be dealt with easily by adding a strong-convexity term to $F_u(\bu)$, thus recovering global convergence guarantee of the 3-block ADMM scheme~\eqref{eq:admm_method}. Therefore, we can replace $F_u(\bu)$ in \eqref{eq:kltgv2_reformulate_admm} with the function
$$ \bar{F_u}(\bu;\sigma) = \lambda\,D_{KL}(\bu) + \chi_{\R_+^n}(\bu) + \frac{\sigma}{2}\|\bu\|^2, $$
with $\sigma>0$ strong convexity parameter. In practice, a very small value of $\sigma$ can be used without affecting the numerical results.
For the 3-block ADMM algorithm applied to the problem
\begin{equation}\label{eq:kltgv2_reformulate_admm_str_convex}
    \begin{array}{ll}
    \min\limits_{\bu,\bw,\bz} & \displaystyle \bar{F_u}(\bu;\sigma) + F_w(\bw) + F_z(\bz) \\
    \st     & M_u\,\bu + M_w\,\bw - \bz = 0,
    \end{array}
\end{equation}
the following convergence result holds (see \cite[Theorem~3.1]{cai:2017}).

\begin{theorem}\label{thm:convergence}
Let $\left\{(\bu^{k},\bw^{k},\bz^{k},\bxi^{k})\right\}_k$ be the sequence generated by applying the ADMM scheme \eqref{eq:admm_method} to problem \eqref{eq:kltgv2_reformulate_admm_str_convex}, where $\sigma>0$ is the strong convexity parameter of $\bar{F_u}(\bu;\sigma)$. Moreover, assume $\rho\in\left(0,\frac{6}{17}\sigma\right)$.
Then, the sequence $\left\{(\bu^{k},\bw^{k},\bz^{k},\bxi^{k})\right\}$ converges to a saddle point $(\bu^*,\bw^*,\bz^*,\bxi^*)$ of the Lagrangian function of problem \eqref{eq:kltgv2_reformulate_admm_str_convex}.
\end{theorem}

It is worth noting that, unlike the classical 2-block ADMM scheme, Theorem~\ref{thm:convergence} imposes a bound on the magnitude of the penalty parameter $\rho$. We observed that in practice such a small value for $\rho$ makes the ADMM method stagnate. Luckily enough, the bound can be relaxed, obtaining fast convergent schemes, as in the case of the numerical experiments shown in this paper.

\subsection{Solving the ADMM subproblems}

Now we show how the solution of the ADMM subproblems can be performed with a small computational effort. First, we note that by exploiting the linearity of the equality constraints, ADMM for the solution of problem \eqref{eq:kltgv2_reformulate_admm_str_convex} can be written in a simplified way. By introducing the scaled Lagrange multipliers $ \bmu^k = \frac{1}{\rho}\bxi^k$, the four update steps of ADMM become
\begin{equation}\label{eq:admm_method_mod_U_problem}
    \bu^{k+1} = \displaystyle \mathrm{arg}\min_{\bu\geq 0} \lambda\,D_{KL}(\bu) + \frac{\sigma}{2}\|\bu\|^2 + \frac{\rho}{2}\left\|M_u\,\bu - \bv^k_u\right\|_2^2,
\end{equation}
with $\bv^k_u = \bz^k - M_w\,\bw^k - \bmu^k$,
\begin{equation}\label{eq:admm_method_mod_W_problem}
    \bw^{k+1} = \displaystyle \mathrm{arg}\min\limits_{\bw\in\R^{2n}} \frac{\rho}{2}\left\|M_w\,\bw - \bv^k_w\right\|_2^2,
\end{equation}
with $\bv^k_w = \bz^k - M_u\,\bu^{k+1} - \bmu^k $,
\begin{equation}\label{eq:admm_method_mod_Z_problem}
\bz^{k+1} = \displaystyle \mathrm{arg}\min\limits_{\bz\in\R^{6n}}\alpha_0\,\|\bz_0\|_{2,1|\R^{2n}} + \alpha_1\,\|\bz_1\|_{2,1|\R^{4n}} + \frac{\rho}{2}\left\|\bz - \bv^k_z\right\|_2^2,
\end{equation}
with $\bv^k_z = M_u\,\bu^{k+1} + M_w\,\bw^{k+1} + \bmu^k$, and
\begin{equation}\label{eq:admm_method_mod_mu_update}
    \bmu^{k+1} = \displaystyle \bmu^k +M_u\,\bu^{k+1} + M_w\,\bw^{k+1} - \bz^{k+1}.
\end{equation}

\subsubsection{Subproblem in $\bu$}
Problem \eqref{eq:admm_method_mod_U_problem} is a bound-constrained smooth optimization problem. We solve it by the Quasi-Newton Projection (QNP) method proposed in \cite{landi:2012}.
To illustrate the idea behind QNP, we consider the Hessian of the objective function, say $f(\bu)$, in \eqref{eq:admm_method_mod_U_problem}, which has the form
\begin{multline*}
    \nabla^2 f(\bu) = \nabla^2 D_{KL}(\bu) + \sigma I_n + \rho M-u^\top M_u = \\
    = A^\top G(\bu)A + \sigma I_n + \rho(D_H^\top D_H + D_V^\top D_V),
\end{multline*}
where
$$ G(\bu) = \diag\left(\bg(\bu)\right),\quad [ \bg(\bu) ]_j = \frac{b_j}{[A\bu+\bgamma]_j^2}. $$
At each step of the QNP algorithm, the Hessian is approximated by the matrix
\begin{equation*}
    \widetilde{H} = \tau(\bu) A^\top A + \sigma I_n + \rho(D_H^\top D_H + D_V^\top D_V),
\end{equation*}
where $\tau(\bu) = \mathrm{mean}_i\left\{[\bg(\bu)]_i\right\}$. Since $A$, $D_H$ and $D_V$ are BCCB matrices, $\widetilde{H}$ is a BCCB matrix too, therefore it can be easily diagonalized by the Discrete Fourier Transform operator $\mF\in\C^{n\times n}$. This allows us to solve the Quasi-Newton system exactly at each step of the QNP algorithm with a small computational cost.

\subsubsection{Subproblem in $\bw$}
Subproblem \eqref{eq:admm_method_mod_W_problem} is an overdetermined least squares problem, where $M_w$ is a full-rank matrix. Hence, we can find its solution by solving the normal equations
\begin{equation}\label{eq:W_subp_norm_eq}
    M_w^\top M_w\,\bw = M_w^\top\bv_w^k.
\end{equation}
$M_w^\top M_w$ can be written blockwise as
\begin{multline}
    M_w^\top M_w = I_{2n} + \mE^\top\mE = \left[\begin{array}{cc} E_{11} & E_{12}\\ E_{21} & E_{22} \end{array} \right] =\\
     {\small = \left[\begin{array}{cc} I_{n} + D_H^\top D_H + \frac{1}{2}D_V^\top D_V & \frac{1}{2}D_V^\top D_H\\[5pt]
        \frac{1}{2}D_H^\top D_V & I_{n} + \frac{1}{2}D_H^\top D_H + D_V^\top D_V \end{array} \right]},
\end{multline}
with $E_{ij}\in\R^{n\times n}$. We note that the four blocks of the matrix have a BCCB structure. Hence, by premultiplying both sides of \eqref{eq:W_subp_norm_eq} by $\begin{bmatrix} \mF & 0\\ 0 & \mF \end{bmatrix}$, we can write it equivalently as
\begin{equation}\label{eq:W_subp_FFTs}
    \left[\begin{array}{cc} \Psi_{11} & \Psi_{12}\\
        \Psi_{21} & \Psi_{22}
    \end{array} \right] \left[\begin{array}{c}
        \mF\,\bw_1\\ \mF\,\bw_2
    \end{array} \right] =
    \left[\begin{array}{c}
        \mF\,\bs_1^k\\ \mF\,\bs_2^k
    \end{array} \right],
\end{equation}
where $\Psi_{ij} = \mF E_{ij}\mF^* \in \C^{n\times n}$ is diagonal, $\bs^k = M_w\bv^k_w$, $\bw = [\bw_1^\top,\,\bw_2^\top]^\top$, and $\bs^k = [(\bs_1^k)^\top,\,(\bs_2^k)^\top]^\top$, with $\bw_1,\bw_2,\bs_1^k,\bs_2^k\in\R^n$. Since $\Psi_{11}$ and $\Psi_{22}$ are invertible, we can use a block inversion formula to compute the inverse of the matrix in \eqref{eq:W_subp_FFTs}, i.e., 
\begin{equation}\label{eq:W_subp_FFTs_inverse}
        \left[\begin{array}{cc}\Psi_{11} &\Psi_{12} \\\Psi_{21} &\Psi_{22} \end{array}\right]^{-1} = \left[\begin{array}{cc}\Omega_1 &-\Omega_1\Psi_{12} \Psi_{22} ^{-1}\\-\Omega_2\Psi_{21} \Psi_{11} ^{-1}&\Omega_2 \end{array}\right],
\end{equation}
where we set 
\begin{equation*}
\begin{split}
    \Omega_1 &= \left(\Psi_{11} -\Psi_{12} \Psi_{22} ^{-1}\Psi_{21} \right)^{-1},\\
    \Omega_2 &= \left(\Psi_{22} -\Psi_{21} \Psi_{11} ^{-1}\Psi_{12} \right)^{-1}.
\end{split}
\end{equation*}
All the blocks of the inverse matrix are diagonal matrices that can be computed in linear time only once before ADMM starts. Thus, to update $\bw$ at each step it is sufficient to compute the vector $\bs^k$ and its Fourier transform, perform the multiplication with the matrix in \eqref{eq:W_subp_FFTs_inverse} and apply the inverse Fourier transform to the output.

\subsubsection{Subproblem in $\bz$}
It is straightforward to check that the minimization problem in \eqref{eq:admm_method_mod_Z_problem} can be split in two separate problems, one in $\bz_0$ and one in $\bz_1$, corresponding to the computation of the proximal operators of the functions
$$f_0(\bz_0) = \frac{\alpha_0}{\rho}\|\bz_0\|_{2,1|\R^{2n}} \;\mbox{ and }\;  f_1(\bz_1) = \frac{\alpha_1}{\rho}\|\bz_1\|_{2,1|\R^{4n}},$$
\noindent respectively.
Since the two $(2,1)$-norms correspond to the sum of $2$-norms of vectors in $\R^2$ and $\R^4$, respectively, the computation of their proximal operators can be split in the computation of $n$ proximal operators of functions that are scaled $2$-norms. 
We recall that, given a vector $\bd$ in $\R^s$ and a constant $c>0$, we have \cite[Chapter~6]{beck:2017book}
\begin{equation}\label{eq:proxi_euclidean_norm}
\begin{split}
   \prox_{c\|\cdot\|_2}(\bd) & = \mathrm{arg}\min\limits_{\by\in\R^s} c\|\by\|_2 + \frac{1}{2}\|\by-\bd\|_2^2 =\\
   & = \max\left\lbrace\frac{\|\bd\|-c}{\|\bd\|},\,0 \right\rbrace\bd. 
\end{split}
\end{equation}
Hence, by defining $\bp^k = \left[[\bv_z^k]_1, \ldots, [\bv_z^k]_{2n}\right]^\top$ and $c_0=\frac{\alpha_0}{\rho}$, we can compute the update of $\bz_0$ by setting
$$\left[[\bz_0^{k+1}]_i,\,[\bz_0^{k+1}]_{n+i}\right]^\top = \prox_{c_0\|\cdot\|_2}\left([p^k_i,\,p^k_{n+i}]^\top\right), $$
for each $i\in\{1,\ldots,n\}$.
In a similar way, by defining $\bq^k = \left[[\bv_z^k]_{2n+1},\ldots,[\bv_z^k]_{6n}\right]^\top$ and $c_1=\frac{\alpha_1}{\rho}$, we can compute the update of $\bz_1$ by setting
\begin{multline*}
\left[[\bz_1^{k+1}]_i,\,[\bz_1^{k+1}]_{n+i},\,[\bz_1^{k+1}]_{2n+i},\,[\bz_1^{k+1}]_{3n+i}\right]^\top =\\  =\prox_{c_1\|\cdot\|_2}\left([q^k_i,\,q^k_{n+i},\,q^k_{2n+i},\,q^k_{3n+i}]^\top\right),    
\end{multline*}
for each $i\in\{1,\ldots,n\}$.

\section{Numerical experiments \label{sec:numerical_experiments}}
\graphicspath{{figures/}}

\begin{figure}[h]
\begin{center}
			\includegraphics[width=4.5cm]{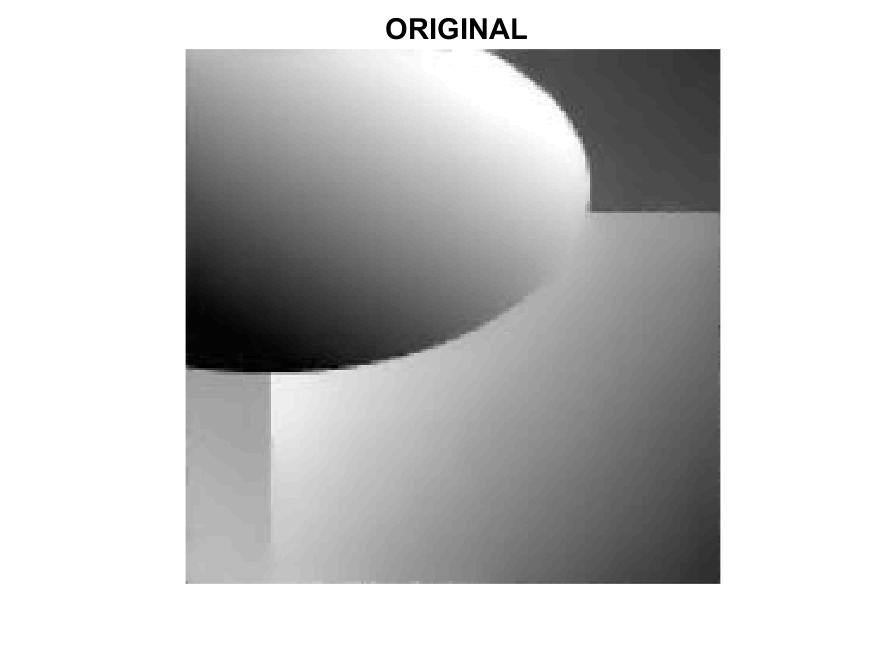} 
			\\[-3mm]
			\includegraphics[width=4.5cm]{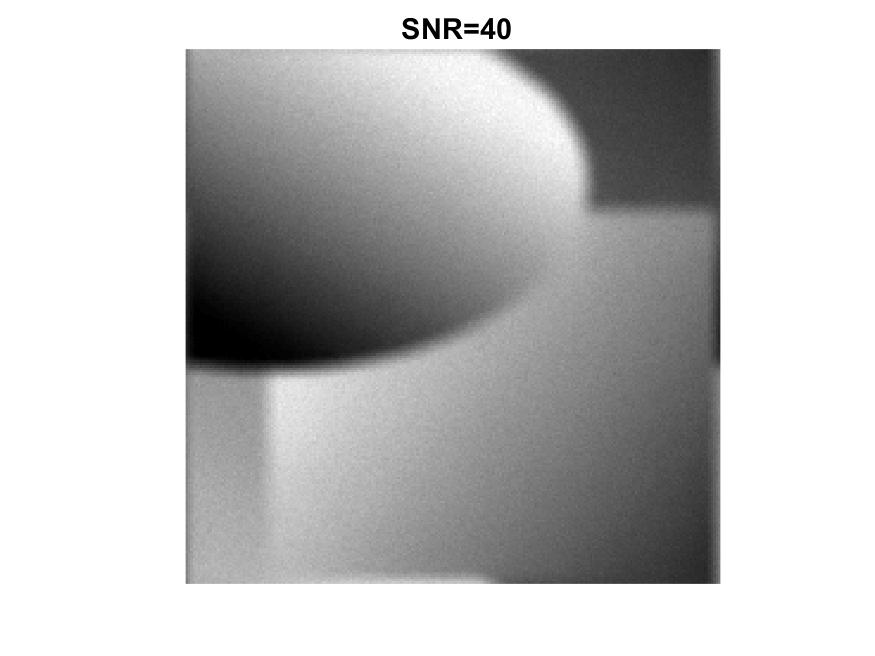}\hspace{-.6cm}
            \includegraphics[width=4.5cm]{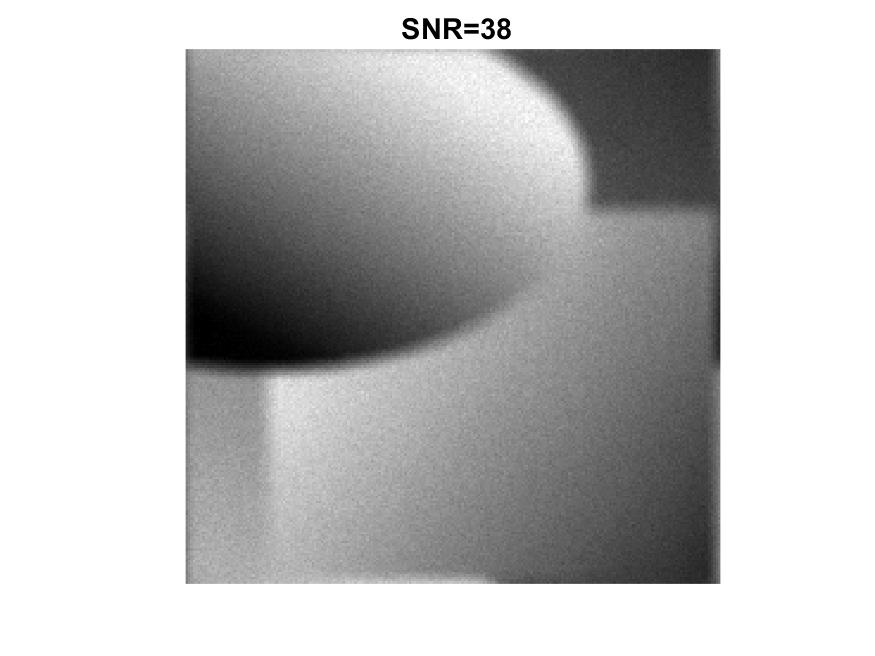}\\[-3mm]
            \includegraphics[width=4.5cm]{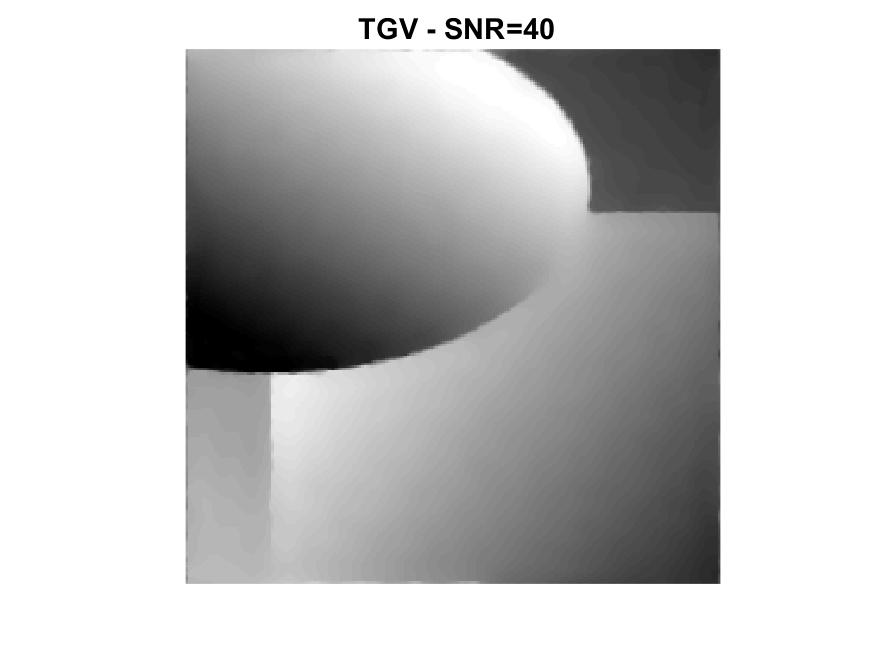}\hspace{-.6cm}
            \includegraphics[width=4.5cm]{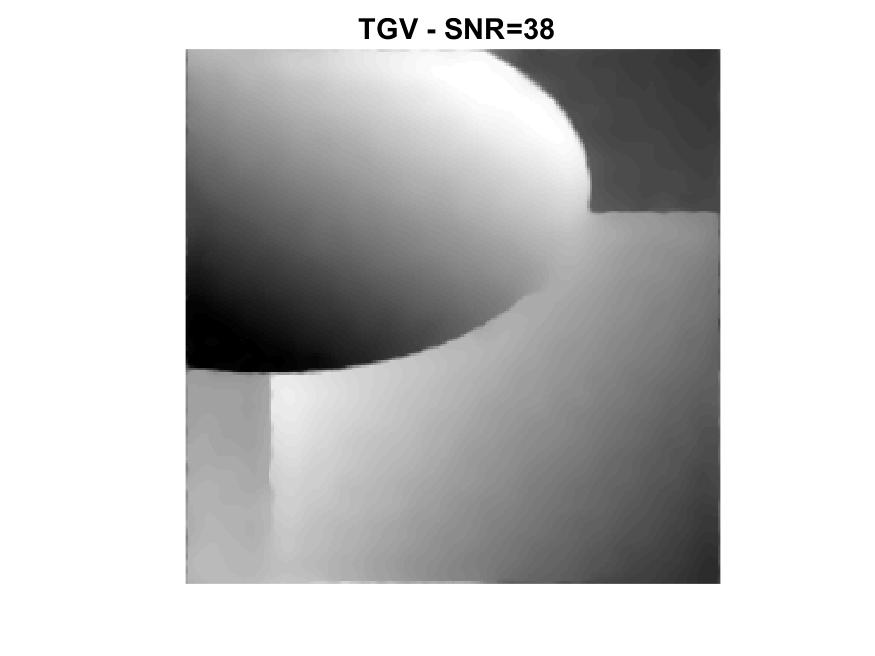}\\[-3mm] 
            \includegraphics[width=4.5cm]{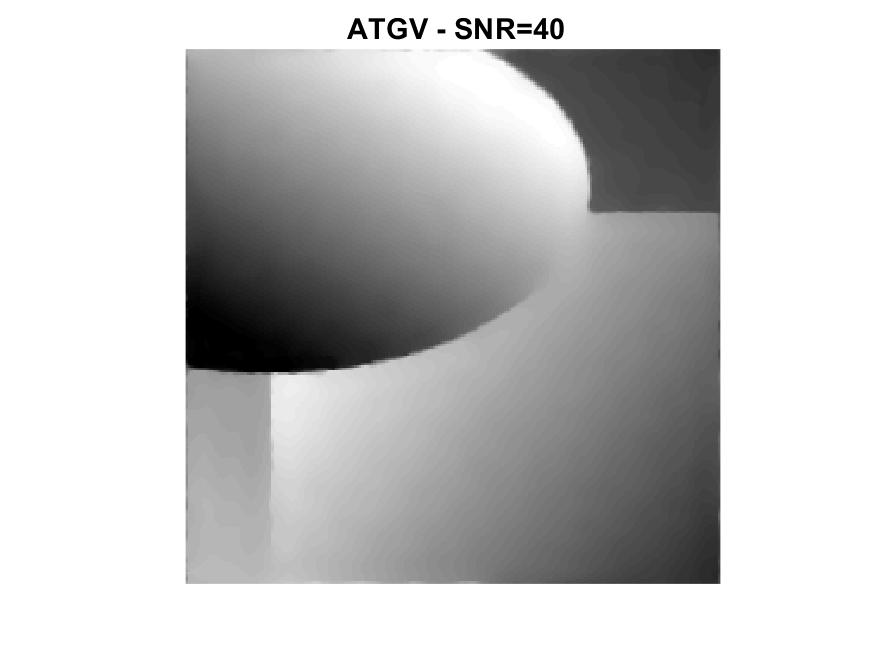}\hspace{-.6cm}
            \includegraphics[width=4.5cm]{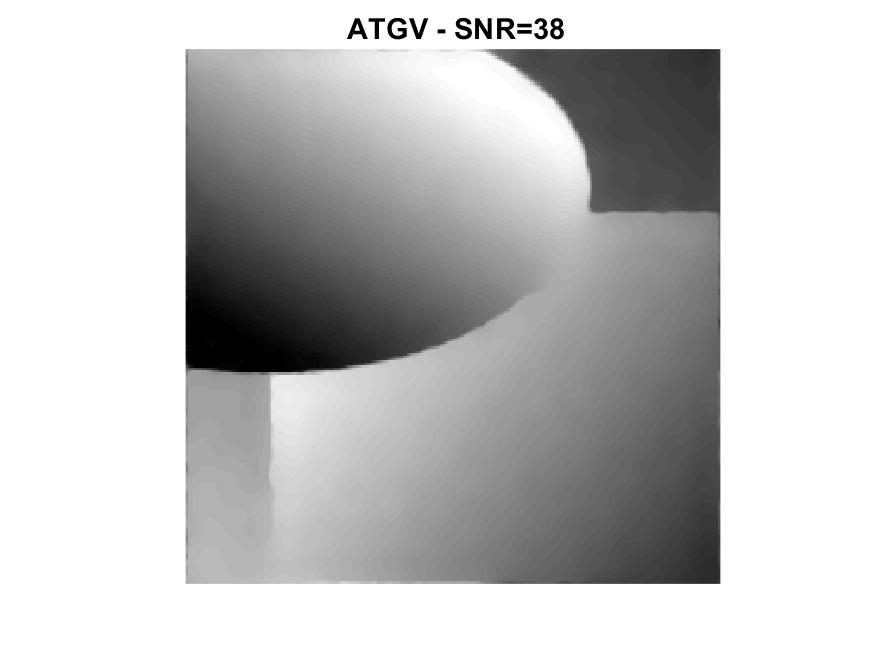}\\[-3mm] 
\end{center}
\caption{Test problem \texttt{phantom} with SNR $=40, 38$.\label{fig:phantom}}
\end{figure}

\begin{figure}[h]
\begin{center}
			\includegraphics[width=5cm]{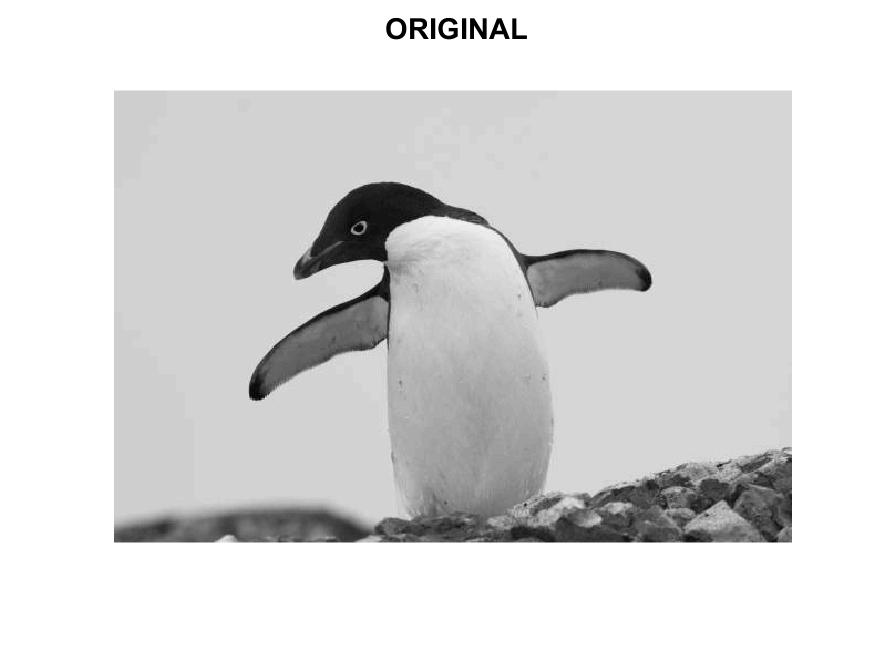} 
			\\[-3mm]
			\includegraphics[width=4.5cm]{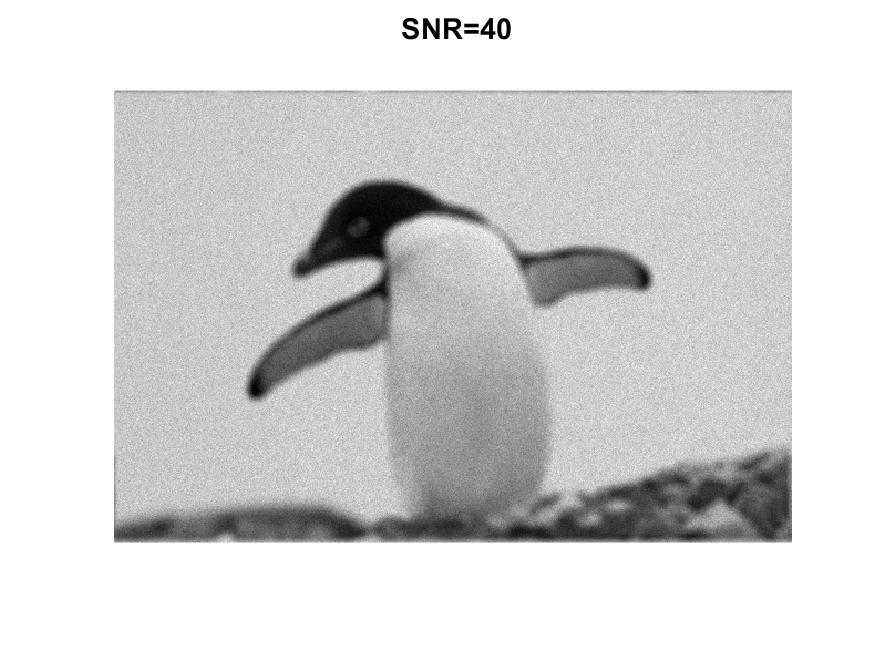}\hspace{-.6cm}
            \includegraphics[width=4.5cm]{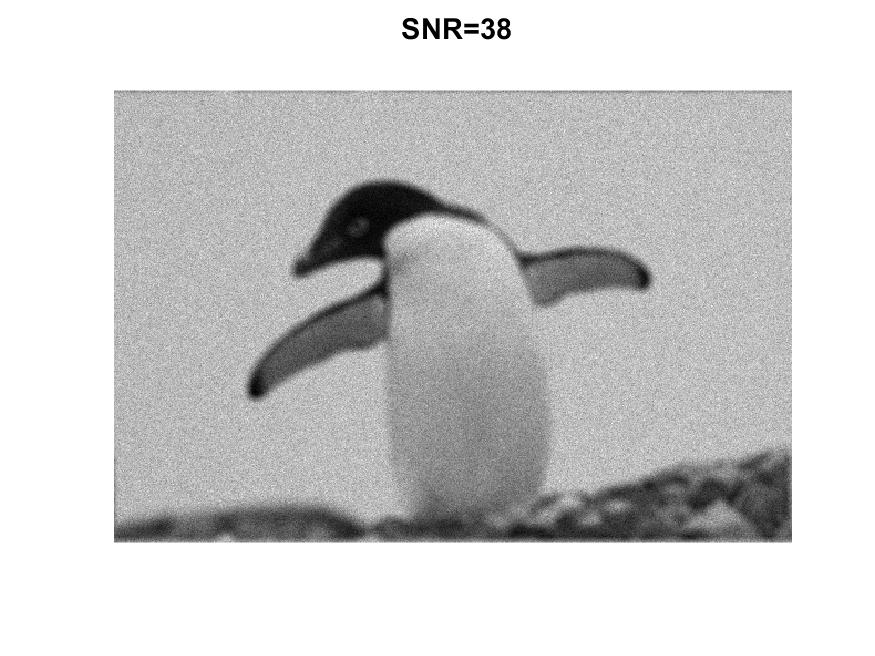}\\[-3mm]
            \includegraphics[width=4.5cm]{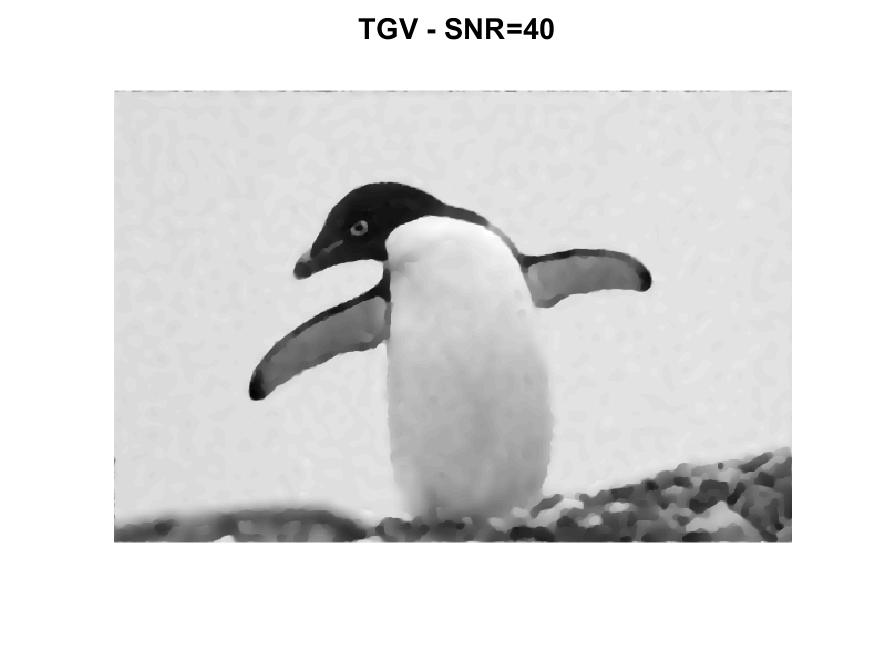}\hspace{-.6cm}
            \includegraphics[width=4.5cm]{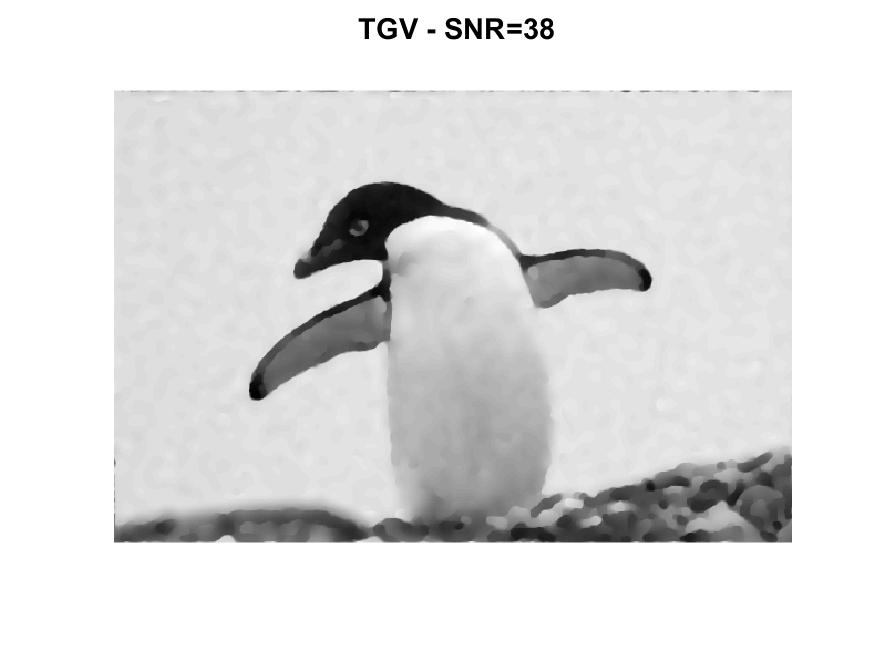}\\[-3mm] 
            \includegraphics[width=4.5cm]{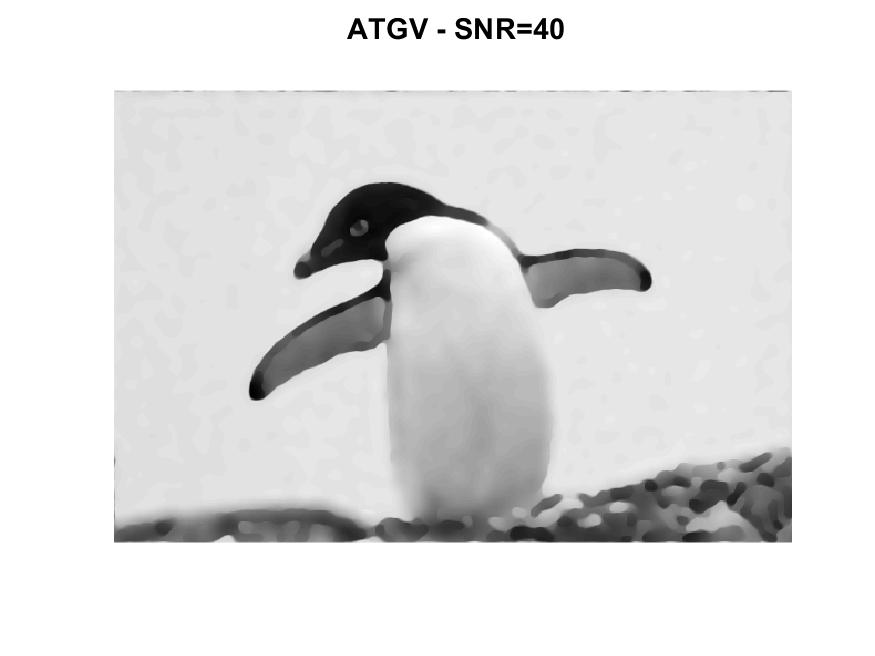}\hspace{-.6cm}
            \includegraphics[width=4.5cm]{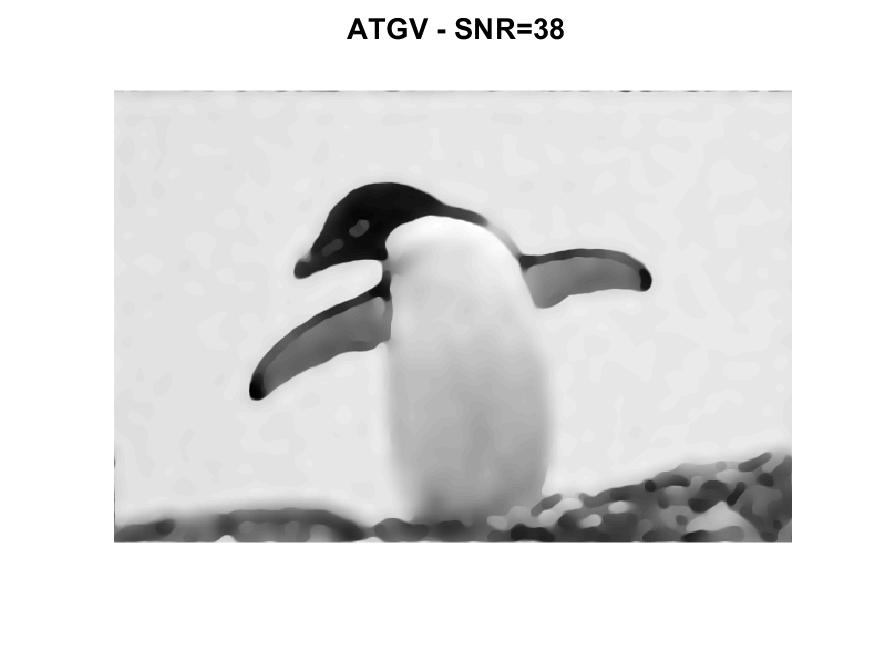}\\[-3mm] 
\end{center}
\caption{Test problem \texttt{penguin} with SNR $=40, 38$.\label{fig:penguin}}
\end{figure}

We compare the results obtained with the 3-block version of ADMM when the value of the regularization parameter is chosen automatically by the proposed strategy and when it is chosen by trial and error. The latter value can be considered as ``almost optimal''.
All the tests were performed using MATLAB R2018a on a 3.50 GHz Intel Xeon E3 with 16 GB RAM and Windows operating system. We use two test images named \texttt{phantom} (size $225\times 225$) and \texttt{penguin} (size $612\times 408$), shown in Figures \ref{fig:phantom} and \ref{fig:penguin}, respectively. These reference images were first blurred by an out-of-focus blur, with radius 5, obtained with the Matlab built-in function \texttt{fspecial}. Then, the Matlab function \texttt{imnoise} was used to add Poisson noise to the blurred images. The intensities of the reference images were pre-scaled to get noisy images with Signal to Noise Ratio (SNR) equal to 42, 40 and 38 dB. For images corrupted by Poisson noise, the SNR is defined as
\begin{equation*}
    \text{SNR} = 10\log_{10}\left(\frac{N_\text{exact}}{\sqrt{N_\text{exact} + N_\text{background}}} \right),
\end{equation*}
\noindent
with $N_\text{exact}$ and $N_\text{background}$ being the total number of photons in the exact image and in the background term, respectively. Finally, the corrupted images were scaled to have the maximum intensity value equal to one. Figures \ref{fig:phantom} and \ref{fig:penguin} also show the blurred and noisy images when SNR $= 40, 38$.

We use the relative error with respect to the original image, computed as $\| \bu^{(k)} - \bu^* \| / \| \bu^* \|$, where $\bu^*$ denotes the original image, to evaluate the quality of the restoration. The time (in seconds) is used as a measure of the computational cost.

The 3-block ADMM described in Section~\ref{sec:3block_admm} was used at each iteration of Algorithm~\ref{alg:fixedpoint} to solve the minimization problem at line~\ref{alg:fixedpoint_line_3}. The stopping condition of ADMM was defined as
$$ \|\bu^{k+1}-\bu^{k}\|/\|\bu^{k+1}\|\leq 10^{-4}. $$
A maximum number of 500 iterations was also fixed. 
The value of the penalty parameter was set as $\rho=1$ and the weights of the TGV regularizer were chosen as $\alpha_0=\beta$ and $\alpha_1=(1-\beta)$ with $\beta=0.1$.

As regards Algorithm \ref{alg:fixedpoint}, the initial guesses for the regularization parameter and the restored image were set as
$$
    \lambda^{0} = 10 \frac{\alpha_0 \left\|\nabla \bb\right\|_{2,1|\R^{2n}}}{D_{KL}(\bb)}
$$
and
$$
    \bu^0 = \bb
$$
(we have implicitly assumed that $\bb$ has been converted into a real vector with entries ranging in the same interval as the entries of $\bu$). The value $\gamma=2.5$ was used. Furthermore, at each iteration $j$, the previous approximation of $\bu$ (namely $\bu^j$) was used as a starting point for the 3-block ADMM. Algorithm~\ref{alg:fixedpoint} was halted when the relative distance between two successive values of $\lambda$ was less than 0.9 or after 5 iterations.

The trial-and-error strategy consisted in running the ADMM algorithm, with initial guess $\bu^0 = \bb$, several times on each test image, varying the value of $\lambda$ at each execution. The stopping criterion for ADMM and the value of $\alpha_0$, $\alpha_1$ and $\rho$ were the same as described above. The value of $\lambda$ corresponding to the smallest relative error at the last iteration was chosen as the ``optimal'' value.

For each test problem, Table \ref{tab:1} reports the values of the regularization parameter obtained with the automatic procedure (ATGV) and the trial-and-error strategy (TGV), the relative error for the corresponding restorations, the number $it_j$ of iterations required by Algorithm~\ref{alg:fixedpoint}, and the execution time in seconds. The numerical results reported in Table \ref{tab:1} show that the automatic procedure produced good quality reconstructions compared with the ``optimal'' ones. Overall, the value of the regularization parameter determined by ATGV appears to be smaller than the trial-and-error one. This, in general, leads to oversmoothing in the restored image. It is worth noting that the oversmoothing of the image determined a better reconstruction quality in the case of \texttt{phantom}, which is piecewise smooth, whereas it led to slightly larger relative errors in the real image \texttt{penguin}. Interestingly, for 1 out of the 6 cases under analysis (\texttt{phantom} with SNR = 42), the initial choice $\lambda_0$ for the regularization parameter yielded a sufficiently good reconstruction in about the same time as the non-automatic procedure. For all the other cases, it is worth noting that the execution time required by the automatic strategy was from 2 to 4.5 times larger than the time needed with the fixed choice of the regularization parameter. We believe this is a quite good performance, considering that a naive trial-and-error strategy may require far more attempts before a reasonable restoration quality is reached. Moreover, in a real-life scenario the real image and the relative error are not available at all, making it impossible to determine $\lambda$ by trial-and-error. 

\begin{figure}[h]
\begin{center}
            \includegraphics[width=4cm]{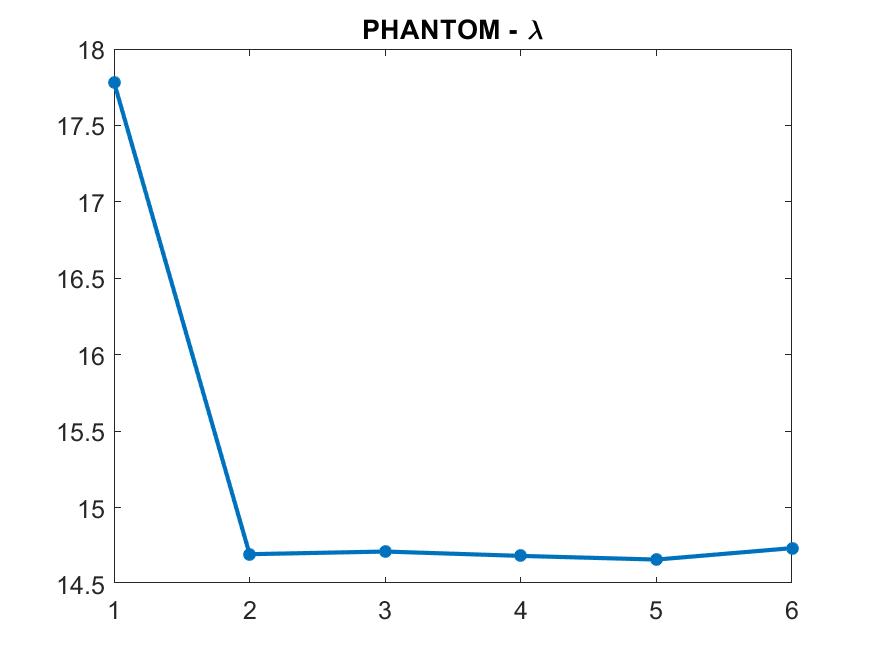}\hspace{-.6cm}
            \includegraphics[width=4cm]{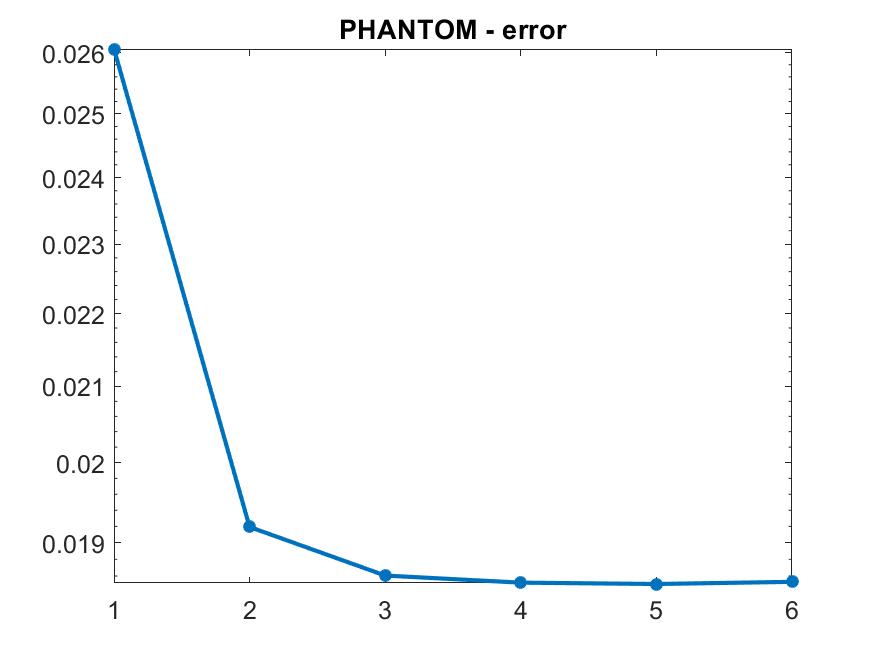}\\
\end{center}
\caption{Test problem \texttt{phantom} with SNR $=38$: values of the regularization parameter (left) and corresponding relative errors (right).\label{fig:phantom2}}
\end{figure}

Figures \ref{fig:phantom} and \ref{fig:penguin} show the restored images when SNR $=40, 38$. Finally, Figure \ref{fig:phantom2} shows the sequence of values of $\lambda$ computed by the automatic procedure and the relative errors of the corresponding images for the \texttt{phantom} test problem with SNR $=38$. It is interesting that the sequence of the regularization parameter values stabilizes after just a couple of iterations.

\begin{table}[h!]
{\scriptsize	
\begin{center}
\begin{tabular}{|c|c|c|c|c|c|}
\hline
SNR  &  Method     & $\lambda$ & Rel. Error  &  $it_j$ &  Time        \\
\hline
\multicolumn{6}{|c|}{\texttt{phantom}} \\
\hline
\multirow{2}{*}{42}  
&  TGV & 31.00 & 1.4057e-02 & -- & 7.78 \\ 
& ATGV & 83.31 & 1.3945e-02 & 1 & 7.21 \\ 
\hline
\multirow{2}{*}{40}  
&  TGV & 20.00 & 1.6447e-02 & -- & 7.72 \\
& ATGV & 17.54 & 1.5940e-02 & 2 & 15.94 \\
\hline
\multirow{2}{*}{38}  
&  TGV & 10.00 & 1.8746e-02 & -- & 7.31 \\
& ATGV & 6.73 & 1.8035e-02 & 2 & 15.16 \\
\hline
\multicolumn{6}{|c|}{\texttt{penguin}} \\
\hline
\multirow{2}{*}{42}  
& TGV & 20.00 & 2.9853e-02 & -- & 32.61 \\ 
& ATGV & 19.65 & 2.9786e-02 & 2 & 73.97 \\ 
\hline
\multirow{2}{*}{40} 
&  TGV & 10.00 & 3.2596e-02 & -- & 43.08 \\
& ATGV & 2.62 & 3.5951e-02 & 3 & 196.82 \\
\hline
\multirow{2}{*}{38} 
&  TGV & 5.00 & 3.6625e-02 & -- & 50.56 \\
& ATGV & 0.95 & 4.0692e-02 & 3 & 221.94 \\
\hline
\end{tabular}
\vspace*{3mm}
\caption{Numerical results for the test problems.\label{tab:1}}
\end{center}
}
\end{table}

\section{Conclusion}
In this work we focused on the restoration of images corrupted by Poisson noise by means of the KL-TGV$^2$ model, which combines the Kullback-Leibler divergence as data fidelity term and the second-order Total Generalized Variation (TGV$^2$) as regularization term. We presented an automatic strategy for the determination of the regularization parameter and a 3-block version of the ADMM method for the solution of the bound-constrained nonsmooth optimization problem corresponding to the KL-TGV$^2$ model. Preliminary numerical experiments show that the automatic regularization strategy, compared with a standard trial-and-error strategy, produces good results in terms of restoration quality in a reasonable computational time.

Future work directions include a theoretical analysis of the automatic strategy, possibly aimed at improving its practical performance. Moreover, an extension to the case of pixelwise regularization parameters will be considered.


\IEEEtriggeratref{18}



%
%

\bibliographystyle{IEEEtran}
\bibliography{IEEEabrv,biblio_kltgv.bib}

\end{document}